\documentclass[12pt,leqno]{article}

\usepackage{amssymb, amsmath, graphicx}

\numberwithin{equation}{section}

\newtheorem{conv}{Convention}[section]
\newtheorem{const}{Construction}[section]

\newtheorem{defn}{Definition}[section]

\newtheorem{obse}{Observation}[section]
\newtheorem{prop}{Proposition}[section]
\newtheorem{rema}{Remark}[section]

\newtheorem{theo}{Theorem}[section]

\newcommand{\cE}{{\cal E}}

\newcommand{\cM}{{\cal M}}
\newcommand{\cN}{{\cal N}}

\newcommand{\cQ}{{\cal Q}}

\newcommand{\cS}{{\cal S}}

\newcommand{\cU}{{\cal U}}

\newcommand{\bbbA}{{\mathbb A}}
\newcommand{\bbbE}{{\mathbb E}}
\newcommand{\bbbG}{{\mathbb G}}
\newcommand{\bbbH}{{\mathbb H}}
\newcommand{\bbbI}{{\mathbb I}}
\newcommand{\bbbL}{{\mathbb L}}

\newcommand{\bbbR}{{\mathbb R}}
\newcommand{\bbbS}{{\mathbb S}}
\newcommand{\bbbT}{{\mathbb T}}
\newcommand{\bbbU}{{\mathbb U}}

\newcommand{\bbbW}{{\mathbb W}}
\newcommand{\bbbZ}{{\mathbb Z}}

\newcommand{\nbhd}{{\rm nbhd}}

\newcommand{\Perm}{{\rm Perm}}

\newcommand{\Stab}{{\rm Stab}}

\newcommand{\Sym}{{\rm Sym}}

\newcommand{\1}{{[\! [}}
\newcommand{\2}{{]\! ]}}

\title{Periodic Euclidean Graphs on Integer Points [Revised]}
\author{ Gregory McColm \\
         Department of Mathematics \& Statistics \\
         University of South Florida \\
         Tampa, FL 33620\\
         (813) 974-9550, fax (813) 974-2700\\
         {\tt mccolm@usf.edu}\\
         URL {\tt http://www.shell.cas.usf.edu/\~{}mccolm}}

\begin{document}

\maketitle

% end of title page

\begin{abstract}
A uniformly discrete Euclidean graph is a graph embedded in a Euclidean space so that there is a minimum 
 distance between distinct vertices.
If such a graph embedded in a $d$-dimensional space is preserved under $n$ linearly independent 
 translations, it is ``$d$-periodic'' in the sense that the quotient group of its symmetry group 
 divided by the translational subgroup of its symmetry group is finite.
We present a refinement of a theorem of Bieberbach:  given a $d$-periodic uniformly discrete 
 Euclidean graph embedded in a $d$-dimensional Euclidean space of symmetry group $\bbbS$, there 
 is another $d$-periodic uniformly discrete Euclidean graph embedded in the same space whose 
 vertices are integer points (possibly modulo an affine transformation) and whose symmetry group has 
 a (not necessarily proper) subgroup isomorphic to $\bbbS$.
\end{abstract}

\medskip

\noindent{\bf Keywords.}
Crystal nets, crystallographic groups, Euclidean graphs of integer points, 
 periodic graphs, symmetry groups of geometric graphs.

\medskip

\noindent{\bf ZBL classification codes.}
05C25, 05C30, 05C62, 05C85, 20E99, 68U05

\vfil

\newpage

\section{Introduction}

Given a Euclidean graph\footnote{Also known as a {\em geometric graph} or an {\em embedded graph} or even an 
 {\em embedded net}.} -- a graph embedded in a Euclidean space -- that is ``periodic'' with respect to some 
 basis of that space, and assuming that it is ``uniformly discrete (i.e., there is a minimum distance between 
 distinct vertices), we will demonstrate the existence of an isomorphic Euclidean graph of at least ``comparable'' 
 symmetry whose vertices are all points on a geometric lattice.
But first, we will define our terms.

We are working with Euclidean graphs.

\begin{defn}\label{firstdef}
A {\rm Euclidean graph} is a graph whose vertices are points in some Euclidean space and whose edges are line 
 segments joining those vertices.
\end{defn}

There are several definitions of periodicity.
One of the more popular definitions (\cite{klee}, see \cite{JKM}) requires a prior definition of the graph's 
 ``symmetry group.''
Recall that a Euclidean space has an isometry group and that the graph itself has an automorphism group.

\begin{defn}
Given a Euclidean graph of automorphism group $\bbbA$ embedded in a Euclidean space of isometry group $\bbbI$, 
 the {\rm symmetry group} of that graph is the group of isometries in $\bbbI$ that induce automorphisms in 
 $\bbbA$.
\end{defn}

There is a tendency to treat the symmetry group and the corresponding group of automorphisms interchangeably.

Returning to periodicity for the graphs we will be considering in this paper, the popular definition may be 
 boiled down to:

\begin{defn}
A Euclidean graph in a $d$-dimensional Euclidean space of translation group $\bbbT$ is {\rm $d$-periodic} if it's symmetry 
 group $\bbbS$ admits a subgroup $\bbbS \cap \bbbT$ generated by $d$ linearly independent translations so that 
 $\bbbS / (\bbbS \cap \bbbT)$ is finite, and if its vertices are all of finite degree.
\end{defn}

Notice that a $d$-periodic Euclidean graph embedded in a $d$-dimensional Euclidean space admits $d$ linearly 
 independent translations that preserve the graph, and thus, from any point in the space, these translations form 
 a {\em geometric lattice} 
 that has the effect of chopping the space in which the graph is embedded into $d$-dimensional parallelopipeds, 
 ``unit cells,'' all of which contain congruent portions of the graph.

We now state the main theorem of this article.

\begin{theo}\label{main}
Suppose that a uniformly discrete Euclidean graph $\cN$ is embedded in a Euclidean space of dimension $d$.
Suppose that it is $d$-periodic and has a symmetry group $\bbbS$.

Then there is a Euclidean graph $\cN'$ embedded in the same space such that:
\begin{itemize}
 \item
  $\cN'$ is isomorphic to $\cN$.
 \item
  If $V'$ is the vertex set of $\cN'$, then $V'$ is a subset of some geometric lattice.
  Viz., there exists an affine transformation $f$ such that $f[V']$ consists of integer points.
 \item
  The symmetry group of $\cN'$ admits a subgroup isomorphic to $\bbbS$.
\end{itemize}
\end{theo}

We will gild the lily somewhat by showing that if no two edges of $\cN$ (taking the edges here as line 
 segments) intersect, then our construction permits the same to be true of the edges of $\cN'$.

Our goal is to set up a description of the given Euclidean graph $\cN$ based on cyclic paths across labeled 
 vertices, and then using the language of words consisting of these labels to compose a system of 
 simultaneous equations;  the solutions to these systems will define a class of Euclidean graphs -- including 
 graphs $\cN'$ with the desired properties.
Here is an outline of this paper.
\begin{itemize}
 \item
  In Section~\ref{NT}, we will set up the algebraic machinery for describing walks through periodic graphs:  
	 these will consist of formal languages and assignments of isometries to words from these formal languages.
 \item
  In Section~\ref{NG}, we will use these languages set up the systems of homogeneous linear equations, taking a finite 
	 set of words and converting them, word by word, into equations.
	We will demonstrate that solutions to those systems exist, and hence that the desired Euclidean graphs exist.
\end{itemize}
This paper is a companion of \cite{ZC}, which introduces the underlying algorithm on which this paper was 
 based.
A preliminary announcement and description of the algorithm appears in \cite{mainpaper}, and some of its theoretical 
 behavior is described in \cite{primary}.

Since this paper lies in the intersection of several quite different fields, there is a variety of 
 extant notations and nomenclatures.
(See, e.g., \cite{dictionary} for a ``dictionary'' for translating back and forth between the nomenclatures of 
  crystallography and graph theory.)
We fix the conventions of this note as follows.
\begin{itemize}
 \item
  Given groups $\bbbG, \bbbH$, let ``$\bbbG \leq \bbbH$'' mean that $\bbbG$ is a subgroup of $\bbbH$.
  If $X$ is a set, let $\Perm(X)$ be the permutation group on $X$;  if $\bbbG \leq \Perm(X)$, then $\bbbG$ 
   is a {\em group of permutations of $X$}.
  \begin{itemize}
   \item
    If $\bbbG$ is a group of permutations of $X$, and if $f, g \in \bbbG$, denote the composition of $f$ 
     and $g$ by $fg$ so that $(fg)(x) = f(g(x))$ for each $x \in X$.
    Given $f, g \in \bbbG$, the {\em left conjugate} of $g$ under $f$ is ${}^f g = f g f^{-1}$;  for a subgroup $\bbbH 
     \leq \bbbG$, let ${}^f \bbbH$ be the subgroup $\{ {}^f h {:} \; h \in \bbbH \}$.
   \item
    If $\bbbG$ is a group of permutations of $X$, and if $Y \subseteq X$ and $g \in \bbbG$, let $g_Y$ be the restriction 
     of $g$ to $Y$;  let $g[Y] = \{ g(y) {:} \; y \in Y \}$.
    We will restrict $\bbbG$ to $Y \subseteq X$ to get a group acting on $Y$ as follows:
      \[
     \bbbG_Y = \{ g_Y {:} \; g \in \bbbG \; \& \; g[Y] = Y \};
      \]
     note that $\bbbG_Y$ is indeed a group of permutations of $Y$.
   \item
    If $\bbbG$ is a group of permutations of $X$ and, for each $x, y \in X$, there exists $g \in \bbbG$ such that 
     $g(x) = y$, we say that $\bbbG$ is {\em transitive} on $X$.
    On the other hand, if $\bbbG$ is not transitive on $X$, and $x \in X$, then the {\em orbit} 
     of $x$ under $\bbbG$ is $\bbbG(x) = \{ g(x) {:} \; g \in \bbbG \}$.
    If $x$ and $y$ are in the same orbit, write $x \sim y$.
   \item
    And if $\bbbG$ is a group of permutations of $X$, and $x \in X$, the {\em stabilizer} of $x$ in $\bbbG$ is 
     $\Stab(\bbbG,x) = \{ g \in \bbbG {:} \; g(x) = x \}$.
		If the stabilizer of $x$ is trivial, we say that $x$ is {\em free} in $\bbbG$.
  \end{itemize}
 \item
  For sets $S$, $T$, let $S - T$ be the set difference $\{ s \in S {:} \; s \not\in T \}$.
  \begin{itemize}
   \item
    If $\cN = \langle V, E \rangle$ is a graph of vertices $V$ and edges $E$, and if $v \in V$, let $\nbhd(v) = 
     \{ v \} \cup \{ w \in V {:} \; \{ v, w \} \in E \}$ be the {\em neighborhood} of $v$.
		If $S \subseteq V$, let 
      \[
     \partial S = \left(\bigcup_{v \in S} \nbhd(v) \right) - S
      \]
     be the {\em boundary} of $S$, so that $\nbhd(v) = \partial\{v\} \cup \{ v \}$.
   \item
    If $S \subseteq V$, the {\em induced subgraph} of $\cN$ of vertex set $S$ is $\cN[S] = \langle S, E [S] \rangle$, 
		 where $E [S] = \{ \{ u, v \} \in E {:} \; u, v \in S \}$.
    Going further, let $\cN^{\partial}[S] = \langle S \cup \partial S, E^{\partial}[S] \rangle$, where 
     $E^{\partial}[S] = \{ e \in E {:} \; e \cap S \neq \varnothing \}$:  $\cN^{\partial}[S]$ gives us $\cN[S]$ plus 
		 those edges that connect from $S$ to $\partial S$, plus the vertices incident to those edges - but not those 
		 edges connecting vertices of $\partial S$ to other vertices of $\partial S$.
  \end{itemize}
 \item
  In this article we will presume $d$ to be some fixed finite dimension.
  Let $\bbbR$ be the set of real numbers, so that $\bbbR^d$ is an $d$-dimensional space, which we will (sloppily) 
   treat both as Euclidean $d$-space and as a vector space.
	Thus we will represent a point by a boldface lowercase Roman letter, and a vector that way as well;  given 
	 two points $\mathbf a$ and $\mathbf b$, let $[\mathbf a, \mathbf b]$ represent the line segment 
	 $\{ \mathbf a + t(\mathbf b - \mathbf a) : 0 \leq t \leq 1 \}$.
  \begin{itemize}
   \item
    Recall (from, say, \cite{Yale}) that an {\em isometry} on $\bbbR^d$ is a map $g \in \Perm(\bbbR^d)$ such 
     that for all $\mathbf x, \mathbf y \in \bbbR^d$, $\| \mathbf x - \mathbf y \| = \| g(\mathbf x) - 
     g(\mathbf y) \|$, where $\| \cdot \|$ is the usual Euclidean metric.
    Let $\bbbI_d$ be the group of isometries on $\bbbR^d$, and let $\bbbT_d$ be the group of translations on $\bbbR^d$;  
     recall (from, say, \cite{Yale} again) that $\bbbT_d$ is a normal subgroup of $\bbbI_d$.
   \item
    A set $S \subseteq \bbbR^d$ is {\em uniformly discrete} if there exists $\varepsilon > 0$ such that for all 
     $\mathbf x, \mathbf y \in S$,
      \[
     \mathbf x \neq \mathbf y \implies \| \mathbf x - \mathbf y \| \geq \varepsilon.
      \]
  \end{itemize}
 \item
  Given a group $\bbbG$ of isometries on $\bbbR^d$, a {\em fundamental region} of $\bbbG$ is a simply connected 
   subset of $\bbbR^d$ that intersects each orbit of $\bbbG$ exactly once.
 \item
  Recall (from, say, \cite{Yale}) that each isometry is expressible in the form $\mathbf x \mapsto \mathbf 
   b + M\mathbf x$ for some fixed vector $\mathbf b \in \bbbR^n$ and some fixed orthonormal matrix $M$ (i.e., 
   the columns of $M$ form an orthonormal basis of $\bbbR^n$);  denote this isometry by $\1 \mathbf b, M \2$, 
   and note that the translation by vector $\mathbf b$ is $\1 \mathbf b, I \2$, $I$ being the identity matrix.
  We can call $\mathbf b$ the {\em vector} component of the isometry and $M$ the {\em linear} component.
	We will employ the notation $\mu(\1 \mathbf b, M \2) = M$.
  \begin{itemize}
   \item
    We compose isometries so:
      \[
     \1 \mathbf b, M \2 \1 \mathbf a, N \2 = \1 \mathbf b + M\mathbf a, MN \2,
      \]
     and hence $\1 \mathbf b, M \2^{-1} = \1 -M^{-1} \mathbf b, M^{-1} \2$.
	  (Recall that orthonormal matrices are nonsingular, and that the transpose of a real orthonormal matrix 
		  is its inverse.)
   \item
    Given a group $\bbbG \leq \bbbI_d$, the {\em point group} of $\bbbG$ is the quotient group 
    \begin{eqnarray*}
     \bbbG / (\bbbG \cap \bbbT_d) &\cong &\{ \1 \mathbf 0, M \2 {:} \; 
                               \mbox{\rm for some $\mathbf a \in \bbbR^d$, $\1 \mathbf a, M \2 \in \bbbG$} \} \\
                                &\cong &\{ M \in \bbbR^n {:} \; 
                               \mbox{\rm for some $\mathbf a \in \bbbR^d$, $\1 \mathbf a, M \2 \in \bbbG$} \} \\
																&=& \mu[\bbbG],
    \end{eqnarray*}
     where $\mathbf 0$ is the zero vector.
    In this article, we will follow the chemists and treat a point group of $\bbbG$ as the corresponding group 
		 $\mu[\bbbG]$ of matrices.
  \end{itemize}
 \item
  A {\em lattice group} in $\bbbR^d$ is a group $\bbbL \subseteq \bbbT_d$ such that for some basis 
   $\mathbf l_1, \ldots, \mathbf l_d$ of $\bbbR^d$, $\bbbL$ is generated by the translations $\1 \mathbf l_1, I \2, 
	 \ldots, \1 \mathbf l_d \2$.
  If $\mathbf b \in \bbbR^d$, then the orbit $\bbbL(\mathbf p) = \{ \mathbf b + \mathbf l {:} \; \mathbf l \in \bbbL \}$ 
	 is the {\em geometric lattice} displaced by $\mathbf b$.
 \item
  Let $\Sigma$ be a finite set.
	Treating the elements of $\Sigma$ as if they were symbols, we can represent finite sequences $(s_1, s_2, \ldots, s_n) 
	 \in \Sigma^{<\omega}$ as {\em strings} so: $s_1 s_2 \cdots s_n$.
  Let $\Sigma^*$ be the set of all strings from $\Sigma$, including the {\em empty string} $\square$.
	The {\em length} of a string is  $|s_1 s_2 \cdots s_n| = n$, and $|\square| = 0$.
	Any set $L \subseteq \Sigma^*$ is a {\em formal language} of {\em alphabet} $\Sigma$.
\end{itemize}

\section{Touring Periodic Graphs}\label{NT}

To prove Theorem~\ref{main}, we will employ a procedure for encoding cycles in Euclidean nets.
We are interested in uniformly discrete $d$-periodic Euclidean graphs in a Euclidean space of dimension $d$;  for the rest of this 
 paper, call such a graph a {\em periodic graph}.

\subsection{The Symmetry Group of a Periodic Graph}\label{symmetry}

We review some basic points of mathematical crystallography (see, e.g., \cite{Coxeter}, \cite{Yale}, and \cite{vinberg}).
We also employ a bit of topology:  given a set $S \subseteq \bbbR^d$, let $\mbox{\rm int}(S)$ be the interior of $S$ and let 
 $\mbox{\rm cl}(S)$ be the topological closure of $S$ (for more ``point-set'' topology, see, e.g., \cite{Hocking}, \cite{Henle}, 
 or \cite{Kuratowski}).
\begin{itemize}
 \item
  A group $\bbbG \leq \bbbI_d$ is {\em crystallographic} if it admits a subgroup $\bbbL$ such that $\bbbG / \bbbL$ is 
	 finite and $\bbbL$ is a lattice group.
	Typically, $\bbbL = \bbbG \cap \bbbT_d$.
 \item
  If $\bbbG$ is crystallographic, then it admits a fundamental region\footnote{This {\em Dirichlet domain} of a free point is 
                                                                                described in \cite{vinberg}.}
   $\Omega \subsetneq \bbbR^n$ such that for some polytope $P$, $\mbox{\rm int}(P) \subsetneq \Omega \subsetneq P$.
  Furthermore, the points in the interior of this fundamental region are {\rm free} in the sense that they are not fixed 
   points of any symmetry of $\bbbG$ besides the identity.
\end{itemize}

Given a periodic graph $\cN$, we obtain its group of symmetries.

\begin{defn}
Given a Euclidean graph $\cN$ embedded in a Euclidean space $\bbbR^n$, a {\rm symmetry} of $\cN$ is an isometry of $\bbbR^n$ 
 whose restriction to the vertex set of $\cN$ is an automorphism of $\cN$.
The group of symmetries is the {\rm symmetry group} of $\cN$.
\end{defn}

We will build a ``scaffolding'' of free and lattice points about the crystal graph.
Note that a periodic graph in $\bbbR^d$ must have a symmetry group whose translation subgroup is generated by $d$ linearly 
 independent translations, whose vector components generate a geometric lattice.
We want to associate a crystal graph with a lattice of points.
We do this as follows.

\begin{conv}
Let $\cN$ be a $d$-periodic graph whose symmetry group $\bbbS$ has a translation subgroup $\bbbS \cap \bbbT_d$ generated by 
 translations of vector components $\mathbf l_1, \ldots, \mathbf l_d$.
We can call $\{ \mathbf l_1, \ldots, \mathbf l_d \}$ a {\rm lattice basis} for the {\rm lattice group} 
    \[
   \bbbL = \{ c_1 \mathbf l_1 + \cdots + c_d \mathbf l_d {:} \; c_1, \ldots c_d \in \bbbZ \}. 
    \]
In addition, choose any point $\mathbf b \in \bbbR^d$.
The set of points 
    \[
   \mathbf b + \bbbL = \{ \mathbf b + c_1 \mathbf l_1 + \cdots + c_d \mathbf l_d {:} \; c_1, \ldots c_d \in \bbbZ \}, 
    \]
 is the {\rm geometric lattice centered at $\mathbf b$}.
\end{conv}

Imagine $\mathbf b + \bbbL$ as a ``scaffold'' of lattice points.

Suppose that $\bbbG$ is crystallographic with lattice group $\bbbL$ of basis $\{ \mathbf l_1, \ldots, \mathbf l_d \}$.
For any $\mathbf b \in \bbbR^d$, $\bbbL$ admits a fundamental region 
  \[
 U_{\mathbf b} = \{ \mathbf b + x_1 \mathbf l_1 + \cdots + x_d \mathbf l_d {:} \; x_1, \ldots x_d \in [-1/2, 1/2) \},
  \]
 whose closure is the parallelopiped $\bar U_{\mathbf b}$.
We call $\bar U_{\mathbf b}$ the {\em unit cell }\footnote{This is what crystallographers call a {\em primitive unit cell}.
	                                                   In this paper, all unit cells will be primitive.
																										(Notice that it is not necessarily true that $U_{\mathbf b}$ is the 
																										 union of (closures of) fundamental regions of $\bbbG$, although 
																										 readers familiar with the IUCr crystallographic tables are used to 
																										 calculating choices of centers $\mathbf b$ that result in such nice 
																										 unit cells.)}
 of $\bbbG$ centered at $\mathbf b$.

Now for a bit of lore.

\begin{rema}
Let $\bbbG$ be crystallographic and fixing $\mathbf b$ and a basis $\{ \mathbf l_1, \ldots, \mathbf l_d \}$ for 
 $\bbbL = \bbbG \cap \bbbT_d$, let $U$ be a fundamental region of $\bbbL$.
\begin{itemize}
  \item
	 First of all, observe that $\bbbL$ is a normal subgroup of $\bbbG$.
	 If $\1 \mathbf m, M \2 \in \bbbG$ and $\1 \mathbf l, I \2 \in \bbbL$, then $\1 \mathbf m, M \2 \1 \mathbf l, I \2 
	  \1 \mathbf m, M \2^{-1} = \1 \mathbf m, M \2 \1 \mathbf l, I \2 \1 -M^{-1} \mathbf m, M^{-1} \2 = 
		\1 M\mathbf l, I \2 \in \bbbG \cap \bbbT_d = \bbbL$.
\end{itemize}
For each $M \in \mu[\bbbG]$, there exists exactly one $\mathbf m_M \in \bbbR^d$ such that $\1 \mathbf m_M, M \2 
 (\mathbf 0) \in U$.
\begin{itemize}
 \item
  There exists a vector $\mathbf m_M$ for each $M \in \mu[\bbbG]$.
	Given $\1 \mathbf m, M \2 \in \bbbG$, choose the unique $\mathbf l \in \bbbL$ such that $\1 \mathbf m, M \2 
   (\mathbf 0) - \mathbf l \in U$;  as $\{ U_{\mathbf l} : \mathbf l \in \bbbL \}$ is a partition of $\bbbR^d$, 
	 one such $\mathbf l$ exists.
	Then set $\mathbf m_M = \mathbf m - \mathbf l$ and $\1 \mathbf m_M, M \2 (\mathbf 0) = \mathbf m - \mathbf l 
	 = \1 \mathbf m, M \2 (\mathbf 0) - \mathbf l \in U$.
 \item
  There exists at most one vector $\mathbf m$ such that $\1 \mathbf m, M \2 (\mathbf 0) \in U$.
	If there were distinct $\mathbf m_1$ and $\mathbf m_2$ such that $\1 \mathbf m_1, M \2 (\mathbf 0), 
	 \1 \mathbf m_2, M \2 (\mathbf 0) \in U$, then as
	  \[
	 \1 \mathbf m_2 - \mathbf m_1, I \2 \1 \mathbf m_1, M \2 = \1 \mathbf m_2, M \2,
	  \]
	 we have $\mathbf m_2 - \mathbf m_1 \in \bbbL$, which is impossible as $\mathbf m_1 = \1 \mathbf m_1, M \2 
   (\mathbf 0)$ and $\mathbf m_1 = \1 \mathbf m_1, M \2 (\mathbf 0)$ are both in $U$ and thus $\mathbf m_2 
	 - \mathbf m_1$ cannot be expressed as an integral combination of $\{ \mathbf l_1, \ldots, \mathbf l_n \}$.
\end{itemize}
Let $\cU = \{ h \in \bbbG : h(\mathbf 0) \in U \} = \{ \1 \mathbf m_M, M \2 : M \in \mu[\bbbG] \}$, and observe 
 that for every $g \in \bbbG$, there exists $\tau \in \bbbL$ and $h \in \cU$ such that $g = \tau h$.
That is because $g = \1 \mathbf m, M \2$ admits a unique $\mathbf l \in \bbbL$ such that $\mathbf m - \mathbf l 
 \in U$, so $\tau = \1 \mathbf l, I \2$ and $h = \1 \mathbf m_M, M \2$ works.
\end{rema}

\subsection{Traveling Through the Graph}

Our goal is a system of linear equations, each representing a cycle in the crystal graph.
We will imagine that each vertex of the graph is a sort of train station, each edge is a track, and that 
 the isometries of the underlying Euclidean space are potential trains that can transport travelers down 
 tracks from one vertex to an adjacent vertex.
{\em We do not require that these train isometries be symmetries of the Euclidean graph}, in fact, in 
 general, they are not, for we can use such an isometry to transport a traveler from a vertex to an 
 adjacent vertex of a different orbit (under the symmetry group).
However, the matrix component of the isometry will be from the point group.

But we want to be very particular about which isometries we use as trains:  in fact, we will choose a 
 fragment of the graph, select as train isometries the translations between vertices in this fragment, 
 and then use conjugates of these isometries for touring the rest of the graph.
We first need to identify this fragment of the graph.
We borrow a notion from geometric group theory.

\begin{defn}\label{TransversalSubgraph}
Given a geometric graph $\cN = \langle V, E \rangle$ on $\bbbR^n$ with a symmetry group $\bbbS$, a {\rm 
 fundamental transversal} is a pair $(S, A)$ such that $S \subseteq V$, $S$ intercepts each $\bbbS$-orbit 
 in $V$ at exactly one vertex, $A$ intersects each $\bbbS$-orbit in $E$ at exactly one edge, every edge 
 in $A$ is incident to at least one vertex in $S$, and for any two vertices in $S$, there is a path along 
 edges in $A$ from one of those vertices to the other.

Let $\kappa {:} \; V \to S$ be the function which, given any $\mathbf v \in S$ and $\mathbf w 
 \in \bbbS(\mathbf v)$, gives us $\kappa(\mathbf w) = \mathbf v$.
\end{defn}

Fundamental transversals are often (confusingly) called {\em fundamental regions} of a graph.
Before we go any further, we should observe that every connected Euclidean graph admits a fundamental transversal.

\begin{prop}[Proposition 2.6 of \cite{DD}]\label{DD}
Each Euclidean graph $\cN$ in $\bbbR^n$ admits a fundamental transversal.\footnote{
  A discrete version of the proof in \cite{Meier} works for Euclidean graphs.
  Start with any vertex $v_0 \in V$ and let $S_0 = \{ v_0 \}$.
  For any $n$, if there are any vertices in $V$ not in an orbit intersecting $S_n$, choose a vertex 
   $v$ in an unrepresented orbit and a path from a vertex in $S_n$ to $v$, and let $v'$ be the first 
   vertex on that path in an orbit unrepresented in $S_n$.
  There must be a vertex $v_{n+1} \in \partial S_n$ such that $v_{n+1} \sim v'$, so let $S_{n+1} 
   = S_n \cup \{ v_{n+1} \}$.
  Continue until all orbits are represented.
	See also \cite{JKM}}
\end{prop}

We assign ``train isometries'' as follows.

\begin{defn}\label{transversalsystem}
Suppose that we are given a Euclidean graph $\cN = \langle V, E \rangle$, with symmetry group $\bbbS$, 
 fundamental transversal $(S, C)$, and $\kappa$.
A {\rm transversal system} is a tuple $\cS = \langle S, \partial S, A, B \rangle$, where
  \[
 A = \{ ( \mathbf u, \mathbf v) \in S \times (S \cup \partial S) {:} \; \{ \mathbf u, \mathbf v \} \in E \},
  \]
 and where $B$ is defined as follows.

For each $\mathbf y \in S \cup \partial S$, fix a symmetry $\1 \mathbf m_{\mathbf y}, M_{\mathbf y} \2 \in \bbbS$ 
 so that $\1 \mathbf m_{\mathbf y}, M_{\mathbf y} \2 (\kappa(\mathbf y)) = \mathbf y$ (if $\kappa(\mathbf y) 
 = \mathbf y$, let $\1 \mathbf m_{\mathbf y}, M_{\mathbf y} \2 = \1 \mathbf 0, I \2$).
Let $B = \{ \1 \mathbf m_{\mathbf y}, M_{\mathbf y} \2 {:} \; \mathbf y \in S \times \partial S \}$.

For any $(\mathbf x, \mathbf y) \in (S \times (S \cup \partial S)) \cap A$, we define the {\rm train isometry} 
 from $\mathbf x$ to $\mathbf y$ as follows.
\begin{itemize}
 \item
  If $(\mathbf x, \mathbf y) \in (S \times S) \cap A$, then the {\rm train isometry} from $\mathbf x$ to $\mathbf y$ 
   is the translation $g_{\mathbf x, \mathbf y} = \1 \mathbf y - \mathbf x, I \2$.
 \item
  If $(\mathbf x, \mathbf y) \in (S \times \partial S) \cap A$, then the {\rm train isometry} from 
   $\mathbf x$ to $\mathbf y$ is 
  \begin{eqnarray*}
   g_{\mathbf x, \mathbf y} 
    &=& \1 \mathbf m_{\mathbf y}, M_{\mathbf y} \2 \1 \kappa(\mathbf y) - \mathbf x, I \2 \\
    &=& \1 \mathbf m_{\mathbf y} + M_{\mathbf y}(\kappa(\mathbf y) - \mathbf x), M_{\mathbf y} \2.
  \end{eqnarray*}
\end{itemize}
\end{defn}

Once we have these symmetries, we can walk from the fundamental transversal through the graph as follows.
First, we assume that the traveler started at some fixed $\mathbf a \in S$, on which $\bbbS$ acts freely.

\begin{defn}
Given $\cN = \langle V, E \rangle$ with unit cell $U_{\mathbf p}$, if there exists $\mathbf v \in V \cap U$ 
 on which $\bbbS = \mbox{\rm Sym}(\cN)$ acts freely, let $\mathbf a = \mathbf v$ and $\cN^{\dagger} = \cN$.

Otherwise, choose $\mathbf v \in V \cap U$ and $\mathbf z \in U$ such that the line segment $[\mathbf a, \mathbf v]$ 
 does not intersect any edge of $E$.
Then let $\cN^{\dagger} = \langle V \cup \bbbS(\mathbf a), E \cup \bbbS([\mathbf a, \mathbf v]) \rangle$, 
 and call $\cN^{\dagger}$ a {\em scaffolded net} derived from $\cN$.
\end{defn}

Notice that given $\cN$, there are many scaffolded nets derived from $\cN$, and in fact, there is more than 
 one isomorphism class of scaffolded nets iff there is more than one orbit of vertices of $\cN$ under 
 $\mbox{\rm Sym}(\cN)$.

\begin{const}\label{KeepSymmetry}
For each $\mathbf x \in V$, we construct an isometry $\1 \mathbf m_{\mathbf x}, M_{\mathbf x} \2 \in \bbbS$ 
 such that $\1 \mathbf m_{\mathbf x}, M_{\mathbf x} \2 (\kappa(\mathbf x)) = \mathbf x$, and thus that 
 $\1 \mathbf m_{\mathbf x}, M_{\mathbf x} \2 \1 \kappa(\mathbf x) -\mathbf a, I \2 (\mathbf a) = \mathbf x$.
\end{const}

For $\mathbf x \in S$, let $\1 \mathbf m_{\mathbf x}, M_{\mathbf x} \2 = \1 \mathbf 0, I \2 \in \bbbS$ and 
 note that in this case, $\1 \mathbf 0, I \2 (\kappa(\mathbf x)) = \mathbf x$.
For $\mathbf x \in \partial S$, let $\1 \mathbf m_{\mathbf x}, M_{\mathbf x} \2 \in B$, and note that in this case, 
  $\1 \mathbf m_{\mathbf x}, M_{\mathbf x} \2 (\kappa(\mathbf x)) = \mathbf x$.
From then on, we proceed by (conjugates of) train isometries.

Suppose that the traveler is at $\mathbf x \not\in S \cup \partial S$, having reached there via 
 $\1 \mathbf m_{\mathbf x}, M_{\mathbf x} \2 \in \bbbS$ such that $\1 \mathbf m_{\mathbf x}, M_{\mathbf x} \2 
 (\kappa(\mathbf x)) = \mathbf x$ (and thus $\1 \mathbf m_{\mathbf x}, M_{\mathbf x} \2 \1\kappa(\mathbf x) 
 - \mathbf a, I \2 (\mathbf a) = \mathbf x$).
Let $\{\mathbf x, \mathbf y\} \in E$ and let $\mathbf y' = \1 \mathbf m_{\mathbf x}, M_{\mathbf x} \2^{-1} (\mathbf y)$, 
 and note that $\mathbf y' \in S \cup \partial S$.
Thus either $\mathbf y' = \kappa(\mathbf y)$ or $\mathbf y' \in \partial S$.
Either way, the traveler may employ a train isometry ${}^{\1 \mathbf m_{\mathbf x}, M_{\mathbf x} \2} 
 g_{\kappa(\mathbf x), \mathbf y'}$ to go from $\mathbf x$ to $\mathbf y$ as follows.
Compute ${}^{\1 \mathbf m_{\mathbf x}, M_{\mathbf x} \2} g_{\kappa(\mathbf x), \mathbf y'} (\mathbf x) = 
 \1 \mathbf m_{\mathbf x}, M_{\mathbf x} \2 g_{\kappa(\mathbf x), \mathbf y'} \1 \mathbf m_{\mathbf x}, M_{\mathbf x} \2^{-1} 
 (\mathbf x) = \1 \mathbf m_{\mathbf x}, M_{\mathbf x} \2 g_{\kappa(\mathbf x), \mathbf y'} (\kappa(\mathbf x)) = 
 \1 \mathbf m_{\mathbf x}, M_{\mathbf x} \2 (\mathbf y') = \mathbf y$.
So to go from $\kappa(\mathbf y)$ to $\kappa(\mathbf x)$ to $\mathbf x$ and then to $\mathbf y$, one may employ 
 the composition
\begin{eqnarray*}
 \1 \mathbf m_{\mathbf y}, M_{\mathbf y} \2 
 &=& {}^{\1 \mathbf m_{\mathbf x}, M_{\mathbf x} \2}g_{\kappa(\mathbf x), \mathbf y'}
     \1 \mathbf m_{\mathbf x}, M_{\mathbf x} \2 \1 \kappa(\mathbf x) - \kappa(\mathbf y), I \2 \\
 &=& \1 \mathbf m_{\mathbf x}, M_{\mathbf x} \2 g_{\kappa(\mathbf x), \mathbf y'} 
     \1 \kappa(\mathbf x) - \kappa(\mathbf y), I \2 \\
 &=& \1 \mathbf m_{\mathbf x}, M_{\mathbf x} \2 \1 \mathbf m_{\mathbf y'}, M_{\mathbf y'} \2.
\end{eqnarray*}
Then as $\1 \mathbf m_{\mathbf x}, M_{\mathbf x} \2$ and $\1 \mathbf m_{\mathbf y'}, M_{\mathbf y'} \2$ are in 
 $\bbbS$, so is $\1 \mathbf m_{\mathbf y}, M_{\mathbf y} \2$.
Furthermore, $\1 \mathbf m_{\mathbf y}, M_{\mathbf y} \2 \1 \kappa(\mathbf y) - \mathbf a, I \2 (\mathbf a) = 
 \1 \mathbf m_{\mathbf y}, M_{\mathbf y} \2 (\kappa(\mathbf y)) = \mathbf y$.
Repeating, we generate the symmetries $\1 \mathbf m_{\mathbf x}, M_{\mathbf x} \2$ satisfying the criteria for 
 the construction.

Note that by taking different routes from $\mathbf \kappa(\mathbf x)$ to $\mathbf x$, we may get different symmetries 
 $\1 \mathbf m_{\mathbf x}, M_{\mathbf x} \2$, but they will all be symmetries in $\bbbS$.
And if two such symmetries are $\1 \mathbf m_{\mathbf x}, M_{\mathbf x} \2$ and 
 $\1 \mathbf m_{\mathbf x}', M_{\mathbf x}' \2$, then $\1 \mathbf m_{\mathbf x}, M_{\mathbf x} \2 
 \1 \mathbf m_{\mathbf x}', M_{\mathbf x}' \2^{-1} \in$ Stab$(\bbbS, \mathbf x)$.

\subsection{Encoding Walks}\label{EncodingWalks}

We will now set up a labeling system to represent these walks through crystal graphs.
To minimize misunderstandings, we will emphasize the distinction between syntax (names of things) versus 
 semantics (things named).
We will employ the nomenclature of {\em formal languages}:  given a finite set $\Sigma$, which we 
 may regard as a set of symbols, a {\em word} or {\em string} from $\Sigma$ is a finite (possibly empty) 
 sequence of elements from $\Sigma$.
Let $\square$ be the ``empty word'' of no elements, and let $x_1 x_2 x_3 \cdots x_n$ represent the 
 $n$-element sequence $x_1, x_2, x_3, \ldots, x_n$ if all these elements are in $\Sigma$.
For any positive integer $n$, let $|x_1 x_2 x_3 \cdots x_n| = n$, and let $|\square| = 0$.
Let $\Sigma^*$ be the set of all words from $\Sigma$, and any subset of $\Sigma^*$ is a {\em language} 
 from $\Sigma$.

The first discrete structure will be a syntactic representation of the transversal subgraph of 
 Definition~\ref{TransversalSubgraph}.
It will assign names to vertices and to each ordered pair of names of adjacent vertices (i.e., to each oriented edge), 
 and it will assign a matrix representing the linear component of an isometry sending the first vertex to the second.

If a traveller is to walk through the graph, there should be an itinerary, so we first turn to the question of how to 
 generate this itinerary, and the decoding scheme for translating itineraries into (compositions of) 
 isometries.
An itinerary will be a word from a language representing a list of instructions of the form, ``assuming 
 that you are in the initial orientation within the initial fundamental transversal on a vertex $\mathbf x \in S$, 
 take the isometry $g_{\mathbf x, \mathbf y}$ to vertex $\mathbf y \in S \cup \partial S$.''

The placement of the traveller, as well as the exact point the traveller is at, is critical.
Imagine a train station with gates at all four points of the compass.
Suppose that the itinerary says, ``take the train at the gate to your left.''
What departing train the traveler takes depends on where the traveler was facing when the traveler read that instruction.
Observe that the point of view of the traveller is determined by the matrix component of the isometry that placed the traveler.
This point of view will be determined by a point group of matrices.

We will need an alphabet to represent movements across edges in the fundamental transversal.
But we will also want to traverse lattice vectors.
So as before, let $S = \{ \mathbf x_0, \mathbf x_1, \ldots, \mathbf x_n\}$ and $\partial S = \{\mathbf x_{n+1}, 
 \ldots, \mathbf x_m\}$.
But now, let $\bbbL$ be the translation subgroup of $\bbbS$, and suppose that $\{ \mathbf l_1, \ldots, \mathbf l_d \}$
 is a basis for $\bbbL$.
Let 
  \[
 \bar \partial S = \{ (\mathbf x_k, \mathbf x_k + \iota \mathbf l_i) {:} \; k = 0, \ldots, n; i = 1, \ldots, d; \iota = -1, 1 \}
  \]
 and notice that there may be repetitions of edges;  to keep our notation from metastasizing, we will keep this 
 nomenclature, repetitions and all.
In fact, for each $k = 0, \ldots, n$, $\iota \in \{ -1, 1 \}$ and $i = 1, \ldots, d$, let $\mathbf x_{m+d(2k+(\iota-3)/2)+i} = 
 \mathbf x_k + \iota\mathbf l_i$.
Then for each $i \leq n$, $j > m$, let $g_{\mathbf x_i, \mathbf x_j}$ be the translation $\1 \mathbf x_j - \mathbf x_i, I \2$.
Notice that if $\iota = -1$, $m + d(2k + (\iota - 3)/2) + i = m + d(2k - 2) + i$, while if $\iota = 1$, then 
 $m + d(2k + (\iota - 3)/2) + i = m + d(2k - 1) + i$.

Let
\[
 \bar A = ((S \times (S \cup \partial S)) \cap E) \cup \bar \partial S
\]
 and let
  \[
 \Sigma = \{ (i, j): (\mathbf x_i, \mathbf x_j) \in \bar A \}
  \]
 and for each $(i, j) \in \Sigma$, let $M_j = M_{\mathbf x_j}$.
If $\# S = \{ 0, 1, 2, \ldots, n \}$ and $\# \partial S = \{ n + 1, n + 2, \ldots, m + 2(n + 1)d \}$, again possibly 
 repeating vertices, we have a digraph $\# \cS = \langle \# S \cup \# \partial S, \Sigma \rangle$.
And define $\#\kappa : \{ 0, 1, \ldots, m \} \to \{ 0, 1, \ldots, n \}$ by
  \[
 \#\kappa(i) = \left\{ \begin{array}{lr} i & i \leq n \\
                                         j & i > n \,\&\, \kappa(\mathbf x_i) = \mathbf x_j.
										   \end{array}
							 \right.
  \]
Notice that for $\mathbf x_{m + d(2k - 1) + i} = \mathbf x_k + \mathbf l_i$, $\kappa(\mathbf x_{m + d(2k - 1) + i}) = 
 \kappa(\mathbf x_k + \mathbf l_i) = \mathbf x_k$, and hence $\#\kappa(m + d(2k - 1) + i) = k$.
Similarly, for $\mathbf x_{m + d(2k - 2) + i} = \mathbf x_k - \mathbf l_i$, $\#\kappa(m + d(2k - 2) + i) = k$.

\begin{defn}
Letting $\Sigma$ be our alphabet, and define the {\em language of legal walks from $0 \in \# S$} by the following recursion.
We will denote this language by $W_0$.
\begin{itemize}
 \item
  $\square$ is the {\em trivial walk} from $0 \in \Sigma$ to itself;  let $\square \in W_0$.
 \item
  Given a legal walk from $0$ to $i$, and given $(\#\kappa(i), j) \in \Sigma$, we can add that last symbol to the string:  
	 if $(0, i_0) (\#\kappa(i_0), i_1) (\#\kappa(i_1), i_2) \cdots (\#\kappa(i_k), i_{k+1}) \in W_0$, and if 
	 $(\#\kappa(i_{k+1}), j) \in \Sigma$, then $(0, i_0) (\#\kappa(i_0), i_1) (\#\kappa(i_1), i_2) \cdots (\#\kappa(i_k), i_{k+1})
	 (\#\kappa(i_{k+1}), j) \in W_0$.
\end{itemize}
\end{defn}

We want to fix where the walks go.

\begin{defn}\label{schedule_defn}
Given the above language of legal walks from $0$ for a given transversal system $\cS = \langle S, \partial S, 
 A, B \rangle$, the {\em schedule of $W_0$ for $\cS$} is the function $\gamma : W_0 \to \bbbI$ defined by the following 
 recursion.
\begin{itemize}
 \item
  First, $\gamma(\square) = \1 \mathbf 0, I \2$.
 \item
  For any $\mathbf s (i, j) \in W_0$, let $\gamma(\mathbf s (i, j)) = {}^{\gamma(\mathbf s)}\gamma_{\mathbf x_i, \mathbf x_j} 
	 \gamma(\mathbf s) = \gamma(\mathbf s) g_{\mathbf x_i, \mathbf x_j}$.
\end{itemize}
\end{defn}
	 
Observe that if $\mathbf s = (0, i_0) (\#\kappa(i_0), i_1) \cdots (\#\kappa(i_k), i_{k+1})$, if $M_l$ is the 
 linear component of $g_{\mathbf x_j, \mathbf x_l}$, then the linear component of $\gamma(\mathbf s)$ is 
 $\mu(\mathbf s) = \prod_{j=1}^{k+1} M_{i_j}$;  notice that $\mu(\mathbf s)$ is independent of the vectors 
 $\mathbf x_{i_1}$, ..., $\mathbf x_{i_{k+1}}$.

We want to assign (archetypic) point groups to the orbits of vertices, and then point groups to the individual 
 vertices, so that two vertices of the same orbit will have conjugate point groups (each a conjugate of the 
 archetypic point group).
Assign $\bbbH_i$ to the orbit of vertex $\mathbf x)i$, $i = 0, \ldots, n$, and assign point groups to the individual 
 vertices by the following recursion, using string $\mathbf s = (0, i_0) (\#\kappa(i_0), i_1)(\#\kappa(i_1), i_2)
 \cdots(\#\kappa(i_k), i_{k+1})$:
\begin{itemize}
 \item
  Let $\bbbH_{\square} = \bbbH_{i_0}$.
 \item
  For each $k$, let $\bbbH_{\mathbf s} = {}^{\mu(\mathbf s)}\bbbH_k$.
\end{itemize}
Notice that if $\mathbf s_1$ and $\mathbf s_2$ terminate at the same vertex of orbit $k$, we will need 
 $\mu(\mathbf s_2)^{-1} \mathbf(s_1) \in \bbbH_k$.

Finally, we want to translate using lattice vectors.
Suppose that $\bbbL$ is the lattice group of $\bbbS$ and that $\mathbf L = \{\mathbf l_1, \ldots, \mathbf l_d\}$ is a basis 
 of $\bbbL$.
Install $\mathbf L$ by adding indices $\{ m + 1, \ldots, m + 2(n + 1)d \}$:  if $k \in \{ 0, \ldots, n \}$ 
 and $i \in \{ 1, \ldots, d \}$ let $m + kd + i$ be the index of the vertex $\mathbf x_k + \mathbf l_i$, and let 
 $m + (n + 1 + k)d + i$ be the index of the vertex $\mathbf x_l - \mathbf l_i$.

\section{Net Generation}\label{NG}

Now that we have our itineraries for traversing the net $\cN$, we use these itineraries to generate ensembles of 
 isomorphic nets.
We will show that among these ensembles will be nets of vertices of integer points (modulo affine transformations, if 
 necessary) of maximal symmetry, completing the proof of Theorem~\ref{main}.
We will proceed as follows.
\begin{itemize}
 \item
  {\em First, we will use the language of legal walks to set up systems of simultaneous equations whose solution spaces 
   define Euclidean graphs (and some symmetries of those nets).}
  We will do this by generating syntactic representations of some paths in the net, and reduce each of these 
   representations to simultaneous systems of homogeneous linear equations whose solutions represent nets.
  In order to do this, we will need to represent not only the steps $(i, j) \in \Sigma$ but also lattice vectors.
 \item
  {\em Second, we need to eliminate those solutions whose corresponding nets would have two distinct vertices at the 
   same place, or two different edges intersecting.}
  We will eliminate these solutions by representing the solutions as ensembles of vector spaces, and then delete these 
   ``bad'' vector spaces from the solution space, leaving only ``good'' vector spaces of solutions for nets 
   with no such collisions or intersections.
\end{itemize}
We will then confirm that if such an adjusted solution space is nonempty, then it has lattice point solutions, 
 and we will be done.

And now for a useful fact about point groups.

\begin{theo}
{\bf (From Bieberbach's Third Theorem, see, e.g., \cite[p. 29]{Sch}, \cite[Theorem 7.1]{charlap} or \cite[Theorem 3.2.2]{Wolf}.)}
For each $d$, there are finitely many isomorphism classes of crystallographic space groups for lattices on $\bbbR^d$, and 
 any two isomorphic crystallographic space groups are affine conjugates.
Thus for each $d$, there are finitely many isomorphism classes of crystallographic point groups (treated as groups of matrices), 
 with any two isomorphic crystallographic point groups being linear conjugates.
\end{theo}

Thus there exist finitely many maximal point groups.
In fact, for $d = 3$, there are two maximal point groups, the symmetry group of the octahedron (denoted O$_h$ in the 
 Sch\"onflies notation, m$\bar 3$m in the Hermann-Maugain notation, and $^*$432 in orbifold notation) and the 
 symmetry group of the hexagonal prism, denoted D$_{6h}$ in in the Sch\"onflies notation, $6/$mmm in the Hermann-Maugain 
 notation, and $^*$622 in orbifold notation.

\begin{rema}\label{BasisChange}
For each $d$, one of the maximal point groups is the point group of the integer lattice $\bbbZ^n$, which consists of 
 integer matrices (in fact, of $n \times n$ matrices whose entries are $-1$, $0$, and $1$).
Every other maximal point group is a linear conjugate of a group of matrices for a linear group on $\bbbZ^n$, i.e., 
 a linear conjugate of a group of integer matrices (these integer matrices need not be orthonormal).
For the rest of this article, we fix a maximal crystallographic point group $\bbbH$, which is a linear conjugate of a group $\bbbH^*$ 
 of integer matrices.
Let's fix this linear equivalence:  for the rest of this article, let $F$ be the matrix such that $\bbbH = {}^F\bbbH^* 
 = F\bbbH^* F^{-1} = \{ FMF^{-1} {:} \; M \in \bbbH^* \}$.
Thus $\bbbH = \bbbH^*$ iff $\bbbH$ is a point group of the integer lattice (i.e., the symmetry group of the $n$-cube) iff $F$ 
 is the identity matrix.
\end{rema}

\subsection{Naive Net Generation}\label{naive}

Having fixed the dimension $d$ of the lattice, the maximal point group $\bbbH$, the conjugacy matrix $F$ and the 
 corresponding conjugate group $\bbbH^*$ of integer matrices, we describe the walks.
In the process, we develop a graph whose vertices are words for legal walks and whose arcs are symbols from $\Sigma$ of 
 Subsection~\ref{NG}.
This is roughly what \cite{Chung} and \cite{klee} called a ``quotient graph''; to avoid confusion (as the term ``quotient 
 graph'' means something quite different elsewhere, e.g., in \cite{DD}), we will call this syntactic object a {\em unit 
 diagram}.

One way to think of the unit diagram is to consider a periodic graph of symmetry group $\bbbS$, whose translation subgroup 
 is the lattice group $\bbbL$.
A  {\em fundamental transversal of $\cN$for the lattice group $\bbbL$} is a set of vertices corresponding to the vertices 
 of the unit diagram: in fact, the unit diagram consists of labels for that transversal.
Given a vertex $\mathbf a$ in $\cN$, the labels of the unit graph are words for legal walks that describe how to obtain a 
 fundamental transversal from $\mathbf a$.

Recall that for a symmetry group $\bbbS$ of a geometric graph $\cN$, the stabilizer of a vertex $\mathbf v$ of $\cN$ 
 is the group Stab$(\bbbS, \mathbf x) = \{ g \in \bbbS : g(\mathbf x) = \mathbf x \}$, and the matrix group of 
 Stab$(\bbbS, \mathbf x)$ is $\mu[\mbox{\rm Stab}(\bbbS, \mathbf x)] = \{ M : \exists \mathbf y (\1 \mathbf y, M \2 
 \in \mbox{\rm Stab}(\bbbS, \mathbf x) \}$.

Given a fundamental transversal of a periodic graph $\cN$ for the translation subgroup of its symmetry group, no two 
 vertices in that transversal may share the same matrix group for their stabilizers.
Otherwise, the translation from one to the other would be in $\bbbS$, contradicting the fact that it is a fundamental 
 transversal for the symmetry subgroup of $\bbbS$.

Fix a periodic graph $\cN = \langle V, E \rangle$ of symmetry group $\bbbS$, whose matrix components form a subgroup of 
 $\bbbH$.
Fix a lattice group $\bbbL$ of translations by vectors $c_1 \mathbf l_1 + \cdots + c_d \mathbf l_d$, $c_1, \ldots, c_d \in 
 \bbbZ$, of periodic graph $\cN$ and primitive unit cell $U$.
Fix a transversal system $\cS = \langle S, \partial S, A, B \rangle$ and as usual supposing w.l.o.g. that $\bbbS$ acts 
 freely on $\mathbf a \in S$, suppose that $\pi : \{ 0, 1, \ldots, m \} \to S \cup \partial S$ has $\pi(0) = \mathbf a$ 
 and $k \leq n$ iff $\pi(k) \in S$.

\begin{const}\label{quotientgraph}
Given the periodic graph $\cN$, we construct a {\rm unit diagram} of $\cN$ as follows.
\end{const}

This algorithm is described in detail in \cite{ZC}.

\medskip

\noindent{\bf Construction.}
We first develop the vertices of the unit diagram.
In the following procedure, we obtain vertices by developing legal walks to reach them, and the sequence of sets of itineraries 
 $Q_0 \subseteq Q_1 \subseteq \cdots$ has as its limit a set of itineraries for the vertex set of a quotient graph, while the 
 sets $\partial Q_0 \subseteq \partial Q_1 \subseteq \cdots$ has as its limit a set of itineraries for the vertices of the 
 boundary of that quotient graph.

We will simultaneously collect a matrix from each word, the matrix being the linear component of the isometry used to 
 reach that vertex from the $0$th matrix.
Let $\cM_0 = \{ I \}$, and we will obtain the sets $\cM_1, \cM_2, \ldots$.
We will also develop a function $\mu$ such that for a word $\mathbf t$, $\mu(\mathbf t) \in \bigcup_i \cM_i$ is the 
 linear component of the isometry of $\mathbf t$.

A language $Q$ is {\em closed under prefixes} if $\mathbf {st} \in Q \implies \mathbf s \in Q$.
Each of the iterates $Q_k$ will be closed under prefixes, and so their union will also be closed under prefixes.
We start with $Q_0 = \{ \square \}$: we start at $\mathbf a$.
Given a set of words for legal walks $Q_k \supseteq Q_{k - 1} \supseteq \cdots \supseteq Q_0$, each of which is closed under 
 prefixes, and given $\mu : Q_k \to \cM_k$, we construct $Q_{k + 1} \supseteq Q_k$ as follows.
(Recall that $(i_0, j_0) (i_1, j_1) \cdots (i_k, j_k)$ represents a legal walk iff $i_0 = 0$ and $\#\kappa(j_l) = i_{l+1}$ 
  for $l < k$.
 Also recall from Definition~\ref{schedule_defn} that $\bbbH_i \leq \bbbH$ is the point group assigned to vertex (orbit) $i$.)
\begin{itemize}
 \item
  Choose a word $\mathbf t \in Q_k$ and a symbol $(i, j)$ such that $\mathbf t (i, j)$ represents a legal walk, and such 
	 that for every $\mathbf t' \in Q_k$, if the last symbol of $\mathbf t'$ is some $(i', j')$ where $\#\kappa(j) = 
	 \#\kappa(j')$, then
	  \[
	 {}^{\mu(\mathbf t (i, j))}\bbbH_{\#\kappa(j)} \neq {}^{\mu(\mathbf t')}\bbbH_{\#\kappa(j')}.
	  \]
	Then let $Q_{k+1} = Q_k \cup \{ \mathbf t (i, j) \}$ and $\cM_{k+1} = \cM_k \cup \{ M_{\mathbf t (i, j)} \}$ and 
	 $\mu(\mathbf t (i, j)) = M_{\mathbf t (i, j)}$.
 \item
  If there is no such word $\mathbf t \in Q_k$ and $(i, j)$, let $Q = Q_k$ and let $\partial Q$ be the set 
	 $\{ \mathbf t (i, j) \in W_0 - Q : \mathbf t \in Q \; \& \; (i, j) \in \Sigma\}$.
\end{itemize}
Since all matrices $\mu(\mathbf t)$ are elements of the finite group $\bbbH$, this recursion eventually halts.
Collapse this set into equivalence classes:
  \[
 \mathbf t \simeq \mathbf t' \quad \mbox{\rm iff} \quad g(\mathbf t)(\mathbf a) = g(\mathbf t')(\mathbf a)
  \]
 and let $Q/{\simeq} = \{ [\mathbf t ]_{\simeq} {:} \; \mathbf t \in \bigcup_m Q_m \}$ and $\partial Q/{\simeq} = 
  \{ [\mathbf t]_{\simeq} {:} \; \mathbf t \in \bigcup_m \partial Q_m \}$.
Thus $Q/{\simeq}$ is the set of classes of itineraries to vertices in a connected subnet, with exactly one vertex from 
 each orbit of $\cN$ in the translation subgroup $\Sym(\cN) \cap \bbbT_n$, while $\partial Q/{\simeq}$ consists of the 
 itineraries to adjacent vertices.

Thus a particular vertex in $\cN$ can be assigned a matrix from a collection of possible orientation matrices, but not necessarily 
 a unique one.
Recall $\gamma : W_0 \to \bbbI$ from Section~\ref{EncodingWalks} and recall from Construction~\ref{KeepSymmetry} that if 
 $\gamma(\mathbf t)(\mathbf a) \in \bbbS(\mathbf a)$, then $\gamma(\mathbf t) \in \bbbS$.
Let 
  \[
 E_Q = \{ \{ [\mathbf t]_{\simeq}, [\mathbf t']_{\simeq} \} {:} \; \{ \gamma(\mathbf t), \gamma(\mathbf t') \} \in E \}.
  \]
We now have a {\em unit diagram} $\cQ = \langle Q, E_Q \rangle$ of equivalence classes of itineraries.
$\blacksquare$

\medskip

From this unit diagram, we will generate a system of simultaneous linear equations whose solutions will represent nets 
 isomorphic to $\cN$ (and degenerate nets that we will delete in the next subsection).
Let's construct these linear equations.

\subsection{Constructing the Linear Forms}

First, we need to construct linear forms from the itineraries.

\begin{const}\label{path}
Given a Euclidean net $\cN$, each itinerary is assigned a corresponding {\em linear form} as follows.
\end{const}

Before we start the construction, we first remark that we will use the matrix components of isometries represented by 
 itineraries.
Recall the linear change of basis matrix $F$ from Remark~\ref{BasisChange} so that the matrices represented in the linear 
 forms are integral -- and hence the linear forms themselves have integer coefficients.
(After we have used the linear forms to generate the systems of equations which we then solve, we will have to reverse 
 the change of basis to obtain the desired nets, by returning from subgroups of the group of integer matrices $\bbbH^*$ 
 to the point group $\bbbH = {}^F\bbbH^*$.)
Notice also that $\bbbH^* = {}^{F^{-1}}\bbbH$.

\noindent{\bf Construction.}
The construction is recursive.
\begin{itemize}
 \item
  Here is the basis of the recursion.
  If $(t, t') \in \Sigma$ and $M = \mu((t, t'))$, let
    \[
   \mathbf x_{(t, t')} = \left( \begin{array}{c} x_1^{(t, t')} \\ \vdots \\ x_d^{(t, t')} \end{array} \right),
    \]
   where $x_1^{(t, t')}, \ldots, x_d^{(t, t')}$ are real-valued variables, and let the {\em form} of $(t, t')$ be 
   the expression 
     \[
    \1 \langle x_1^{(t, t')}, \ldots, x_d^{(t, t')}\rangle, {}^{F^{-1}}M \2.
     \]
  Notice that this form is a syntactic object set up to represent an affine map with an integral linear component.
 \item
  Here is the recursive step.
  If an itinerary $\mathbf t$ has form $\1 \mathbf x, {}^{F^{-1}}M \2$ and $(t, t')$ has form $\1 \mathbf x', {}^{F^{-1}}M' \2$, 
   then the form of $\mathbf t (t, t')$ will be $\1 \mathbf x + {}^{F^{-1}}M\mathbf x', {}^{F^{-1}}(MM') \2$, where addition 
   and multiplication of variable vectors is performed in the usual way.
  Notice that the vector component of a form of an itinerary will be an $d$-dimensional vector whose components will be 
   linear combinations involving a large number of variables.
\end{itemize}
Continue the recursion.
Since each step $(t, t')$ is from a fixed alphabet $\Sigma$, the number of variables here is at most $d|\Sigma|$.
$\blacksquare$

\medskip

It will turn out that we only need forms from itineraries of the unit diagram.
Now that we have these linear forms, we can use them to build homogeneous linear equations as follows.

\begin{const}
Given a scaffolded graph $\cN^{\dagger}$ of a crystal graph $\cN$ with itineraries and linear forms as above, and with a 
 symmetry group $\bbbS$, construct a system of linear equations from the itineraries of the quotient graph as whose solutions 
 will be lattice vectors and positions of vertices of scaffolded graphs isomorphic to $\cN^{\dagger}$, including 
 $\cN^{\dagger}$ itself, possibly with {\rm vertex collisions} (i.e., two or more distinct vertices assigned to the same point; 
 we will deal with vertex collisions later).
\end{const}

\noindent
{\bf Construction.}
Let $\mathbf u_1, \ldots, \mathbf u_n$ be vectors of variables $\mathbf u_i = \langle u_{i,1}, \ldots, u_{i,n} \rangle$, 
 where each $u_{i,j}$ is real-valued;  the motivation is that these variable vectors will range over lattice vectors.
We will use $\mathbf x$ to represent variable vectors ranging over positions of transversal vertices, 
 or to represent sums of products of integer matrices with variable vectors.
Recall that $Q$ is a set of equivalence classes of itineraries, two itineraries being equivalent if they start at (the vertex 
 $\mathbf a$ represented by) $t_0$ and terminate at the same vertex.
Let $\mathbf x_0$ be a vector-valued variable that assumes the constant value $\mathbf a$ in $\cN^{\dagger}$.
\begin{itemize}
 \item
  First, we build Construction~\ref{path} into our system of equations.
	For $\mathbf t' = \mathbf t(i, j)$, using $\1 \mathbf x, {}^{F^{-1}}M \2$ for $\mathbf t$ and $\1 \mathbf x_{i,j}, 
	 {}^{F^{-1}}M_{i,j} \2$ for $(i, j)$, if
	  \[
	 \mathbf v_{i,j} + {}^{F^{-1}}M_{i,j} \mathbf v = \mathbf v' + \sum_{i=1}^d c_i \mathbf l_i
	  \]
	 in $\cN$, then for vector-valued variables $\mathbf x'$, we have the equation $\cE_{\mathbf t \to \mathbf t'}$ representing
	  \[
	 \1 \mathbf x, {}^{F^{-1}}M \2 \1 \mathbf x_{i,j}, {}^{F^{-1}}M_{i,j} \2 = 
		   \1 \sum_{i=1}^d c_i \mathbf u_i, I \2 \1 \mathbf x', {}^{F^{-1}}M' \2,
		\]
	 or, more precisely, as $M M_{i,j} = M'$,
    \[
	 \cE_{\mathbf t \to \mathbf t'} = 
	   \mbox{\rm ``}\mathbf x_{i,j} + {}^{F^{-1}}M_{i,j} \mathbf x = \mathbf x' + \sum_{i=1}^d c_i \mathbf u_i\mbox{\rm ''}.
		\]
	 Notice that as the isometries $\1 \mathbf x_{i,j}, {}^{F^{-1}}M_{i,j} \2$ are not necessarily elements of $\bbbS$, this 
	  is not a multiplication table for $\bbbS$;  that comes below.
 \item
  Recall that two itineraries are {\em equivalent} if, starting from $\mathbf a$, they terminate at the same vertex or lattice 
	 point.
  For each $q \in Q$ and each $\mathbf t, \mathbf t' \in q$, of forms $\1 \mathbf x, {}^{F^{-1}}M \2$ and $\1 \mathbf x', 
   {}^{F^{-1}}M' \2$, the {\em equation witnessing the equivalence of $\mathbf t$ and $\mathbf t'$} is 
    \[
   \cE_{\mathbf t \leftrightarrow \mathbf t'} = \mbox{\rm ``}\mathbf x - \mathbf x' = \mathbf 0\mbox{\rm ''}.
    \]
 \item
  For each $\mathbf t$ and $\mathbf t'$ such that $\gamma(\mathbf t)(\mathbf a) - \gamma(\mathbf t')(\mathbf a) \in \bbbL$, 
	 let $\mathbf t$ have form $\1 \mathbf x, {}^{F^{-1}}M \2$ and let $\mathbf t'$ have form $\1 \mathbf x', {}^{F^{-1}}M' \2$.
  Suppose that $\gamma(\mathbf t)(\mathbf a) - \gamma(\mathbf t')(\mathbf a) = c_1 \mathbf l_1, + \cdots + c_n \mathbf l_n$ 
   (and notice that $c_1, \ldots, c_n \in \bbbZ$).
  Then the {\em equation witnessing the lattice displacement of $\mathbf t$ to $\mathbf t'$ by $\gamma(\mathbf t')(\mathbf a) - 
   \gamma(\mathbf t)(\mathbf a)$} is 
    \[
   \cE_{\mathbf t' \mapsto \mathbf t} 
    = \mbox{\rm ``$\mathbf x - \mathbf x' - \sum_{i=1}^n c_i \mathbf u_i = \mathbf 0$''}.
    \]
   where $\mathbf u_1, \ldots, \mathbf u_n$ are $n$-tuples of real variables.
 \item
  In particular, recalling that we are dealing with a scaffolded net, in which $\mathbf a$ is the initial scaffolding 
   point, for any itinerary $\mathbf t$ for a path from $\mathbf a$ to another scaffolding point $\gamma(\mathbf t)
   (\mathbf a)$ in a(nother) unit cell, we have an equation
    \[
   \cE_{\mathbf a \mapsto \mathbf t} 
    = \mbox{\rm ``$\mathbf x - \sum_{i=1}^n c_i \mathbf u_i = \mathbf x_0$''},
    \]
   where $\1 \mathbf x, {}^{F^{-1}}M \2$ is the form of $\mathbf t$, $\1 \mathbf x_0, I \2$ is the form for $\square$ (where 
   $I$ is the identity matrix) and where the integer coefficients $a_1, \ldots, a_n$ 
   are computed from the original scaffolded net $\cN^{\dagger}$.
\end{itemize}
Notice that as $\cN$ and hence $\cN^{\dagger}$ are periodic, it suffices to do the above for all itineraries $\mathbf t$ in 
 which no $(i, j)$ appears more than once.
Thus we have finitely many equations.

We now add some equations to build in the symmetries of $\cN$.
But first, recall from Subsection~\ref{symmetry} that if $\bbbL = \bbbS \cap \bbbT_d$ has basis $\{\mathbf l_1, \ldots, \mathbf l_d\}$, 
 then, if $U$ is a parallelopiped bounded by those vectors, we can generate $\bbbS$ from
  \[
 \cU = \{ \1 \mathbf b, M \2 \in \bbbS : \mbox{\rm $\mathbf b \in U$} \}.
  \]
(Caution:  Recall that $U$ is not necessarily any kind of ``standard'' unit cell of $\bbbS$ that can be partitioned into convex 
  fundamental regions of $\bbbS$, merely a primitive unit cell of the lattice.)
As $\cU$ is a set of representatives from $\bbbS /  \bbbL$, every symmetry of $\cN$ can be expressed as a composition $\tau g$, 
 where $g \in \cU$ and $\tau \in \bbbL$.

We would also like to identify each point of $\bbbR^d$ with corresponding points in the initial unit cell $U$ in the same 
 orbit under the lattice group $\bbbS \cap \bbbT_d$.
(Notice that $\bbbL$ is fixed by $U$ as it is generated from $U$'s edge vectors.)
Let $\mathbf r \mod U$ be the unique $\mathbf r' \in U$ such that $\mathbf r - \mathbf r' \in \bbbL$.

Finally, we encode the multiplication table of $\bbbS$ itself.
By the periodicity of $\cN$ and hence of $\bbbS$, all we need do is encode the multiplication table for $\cU$ modulo translations, 
 as in \cite{Buser}.
(List the elements of $\cU$ as $g = \1 \mathbf y_g, M_g \2$.)
For each $i$, $j$, $k$ such that $g_i g_j = \tau g_k$ for some $g_i, g_j, g_k \in \cU$ and some 
 $\tau = \1 \sum_{l=1}^d c_{i, j, k, l} \mathbf l_i, I \2 \in \bbbL$, we have the equation
  \[
 \1 \mathbf y_{g_i}, {}^{F^{-1}}M_{g_i} \2 \1 \mathbf y_{g_j}, {}^{F^{-1}}M_{g_j} \2 = 
 \1 \sum_{l=1}^d c_{i, j, k, l} \mathbf l_i, I \2 \1 \mathbf y_{g_k}, {}^{F^{-1}}M_{g_k} \2,
  \]
 or, more precisely, as $M_{g_i} M_{g_j} = M_{g_k}$,
  \[
 \cE_{i,j;k} = 
\mbox{\rm ``}\mathbf y_{g_i} + {}^{F^{-1}}M_{g_i} \mathbf y_{g_j} - \mathbf y_{g_k} - \sum_{l=1}^d c_{i, j, k, l} \mathbf u_i 
            = \mathbf 0\mbox{\rm ''},
  \]
 which suffices since we already have, by construction, $M_{g_i} M_{g_j} = M_{g_k}$.

The result is a system $\cE$ of homogeneous linear equations with integer coefficients.
$\blacksquare$

\medskip

\begin{obse}\label{coeffs}
Since the displacements of $\cN^{\dagger}$ itself gives a solution to this system, this system is soluble.
Since the coefficients are integers, there is a solution consisting of rationals, i.e., all vertices 
 $\mathbf v'$ and vector components $\mathbf b'$ are rational combinations of lattice vectors 
 $\mathbf l'$ from the lattice $\bbbL'$ generated by the values for $\mathbf u_1, \ldots, \mathbf u_d$ 
 for $\cN'$.
\end{obse}

We now posit a solution to this system of simultaneous equations.

\begin{defn}
Let the assignments
\begin{eqnarray*}
 \mathbf x_i &\to & \mathbf v_i', \qquad i = 0, \ldots, m, \\
 \mathbf y_{g_j} &\to & \mathbf b_j', \qquad j = 0, \ldots, p \\
 \mathbf u_k &\to & \mathbf l_k', \qquad k = 0, \ldots, n
\end{eqnarray*}
 be a solution to the simultaneous system of equations $\cE$ so that:
\begin{itemize}
 \item
  $\mathbf v'_0$ is free, and $\mathbf v'_0, \ldots, \mathbf v'_{m'}$ are the positions of the vertices of the interior of the 
   quotient graph, while $\mathbf v'_{m'+1}, \ldots, \mathbf v'_m$ are positions of vertices on the boundary, and 
 \item
  $\mathbf l'_1, \ldots, \mathbf l'_n$ are lattice vectors, and let
    \[
   L' = \{ a_1 \mathbf l'_1 + \cdots + a_n \mathbf l'_n {:} \; a_1, \ldots, a_n \in \bbbZ \}.
    \]
   which is preserved by the subgroup $\bbbL' = \langle g_0', \ldots, g_p'\rangle$, where $g_j' = \1 \mathbf b_j', M_{g_j} \2$, and
 \item
  letting 
    \[
   V' = \{ \mathbf v' + \mathbf l' {:} \; \mathbf v' \in V_0' \; \& \; \mathbf l' \in L' \}
    \]
   be the vertices of the resulting net, and letting
    \[
   E' = \{ \{ \mathbf v' + \mathbf l', \mathbf v'' + \mathbf l'\} {:} \;
           \{ \mathbf v', \mathbf v''\} \in E_0' \; \& \; \mathbf l' \in L' \}.
    \]
\end{itemize}
 then the net $\cN' = \langle V', E' \rangle$ is a {\rm solution} to the system $\cE$.
\end{defn}

We will want $\cN'$ to be isomorphic to $\cN^{\dagger}$, and $\bbbS \leq \mbox{\rm Sym}(\cN^{\dagger})$.

\begin{defn}
Let $\eta$ be the mapping from $\cN$ to $\cN'$ as follows.
$\eta : V \to V', E \to E'$ is induced by the equations above:  $\mathbf v_i \mapsto \mathbf v_i'$ as they are the solutions 
 (for $\cN^{\dagger}$ and $\cN'$, respectively) for the variables $\mathbf x_i$, $\mathbf b_g \mapsto \mathbf b_g'$ as they 
 are the respective solutions for the variables $\mathbf y_g$, and $\mathbf l_j \mapsto \mathbf l_j'$ as they are the 
 respective solutions for the variables $\mathbf u_j$.
Extend $\eta$ to a mapping from all of $\cN^{\dagger}$ to $\cN'$ using
  \[
        \1 \sum_{i=1}^d c_i \mathbf l_i, I \2 \1 \mathbf b, {}^{F^{-1}}M \2 (\mathbf v) \mapsto
        \1 \sum_{i=1}^d c_i \mathbf l_i', I \2 \1 \mathbf b', {}^{F^{-1}}M' \2 (\mathbf v'),
	\]
 or, more precisely,
  \[
 \eta : \mathbf b + \sum_{i=1}^d c_i \mathbf l_i + {}^{F^{-1}}M\mathbf v \mapsto
        \mathbf b' + \sum_{i=1}^d c_i \mathbf l_i' + {}^{F^{-1}}M'\mathbf v'.
	\]
Call $\eta$ the {\em mapping induced by the itineraries}.
\end{defn}

Notice that $\eta$ maps vertices to vertices (and, as we shall see, edges to edges), but does not necessarily extend to 
 a nice map from $\bbbR^d$ to $\bbbR^d$.
We have a nuisance to watch for, a nuisance that will play a part in showing that $\cN^{\dagger}$ is isomorphic to $\cN'$ and 
 that $\bbbS \leq \mbox{\rm Sym}(\cN')$.
Recall from Definition~\ref{firstdef} that the edges of a Eulerian graph are the line segments joining adjacent pairs of vertices.
Say that two edges {\em collide} if they intersect at a point that is in the interior of one or both of them.
(Recall that two vertices collide if they are distinct but at the same point.)
We combine these two notions as follows.

\begin{defn}\label{degenerate}
A Euclidean graph $\cN'$ is {\rm degenerate} if it admits four distinct vertices $\mathbf v_{11}, \mathbf v_{12}, 
 \mathbf v_{21}, \mathbf v_{22}$ and edges $[ \mathbf v_{11}, \mathbf v_{12} ]$ and $[\mathbf v_{21}, \mathbf v_{22}]$ 
 such that $[ \mathbf v_{11}, \mathbf v_{12} ] \cap [\mathbf v_{21}, \mathbf v_{22}] \neq \varnothing$.
\end{defn}

We will deal with degeneracy in Subsection~\ref{degeneracy}, but we will need it as a hypothesis here.
Recall that $\bbbS$ is the symmetry group of $\cN^{\dagger}$, and we get:

\begin{prop}\label{main1}
Suppose that $\cN'$ is not degenerate.
Then $\cN'$ is isomorphic to $\cN^{\dagger}$ using the mapping $eta$.
Further, the lattice translation of $\cN^{\dagger}$ and the symmetries in $\bbbU$ correspond to lattice translations of $\cN'$ 
 and symmetries of $\bbbS'$ of the central unit cell of the lattice of $\cN'$.
Thus the symmetry group $\bbbS$ is isomorphic to a (not necessarily proper) subgroup of $\bbbS'$.
\end{prop}

\noindent{\bf Proof.}
First of all, as $\cN'$ is not degenerate, there are no vertex or edge collisions, so $\eta$ is one-to-one.

For the first sentence, we claim three things:
\begin{enumerate}
 \item
  The mapping $\eta$ induced by the itineraries preserves the edges and is thus a graph homomorphism.
 \item
  All edges of $\cN'$ are derived from edges of $\cN$ via $\eta$, which is thus onto.
 \item
  Sym$(\cN^{\dagger})$ can be embedded in Sym$(\cN')$: $\bbbS \leq$ Sym$(\cN')$.
\end{enumerate}
So $\eta$ is an onto, one-to-one graph homomorphism - hence a graph isomorphism - that preserves the symmetry group of 
 $\cN^{\dagger}$.

First, note that by the definition, $\eta$ preserves vertices and edges in mapping transversal onto transversal 
 (and lattice basis to lattice basis).
Recall that given lattice group $\bbbL$, $\mathbf r \mod U$ is the unique $\mathbf r' \in U$ such 
 that $\mathbf r - \mathbf r' \in \bbbL$.
Given itineraries $\mathbf t_1$, $\mathbf t_2$, where $[ \gamma(\mathbf t_1)(\mathbf a), \gamma(\mathbf t_2)(\mathbf a) ] 
 \in E$, for each $i = 1, 2$, applying $\gamma$ to $\cN^{\dagger}$
  \[
 \mathbf v_i = \gamma(\mathbf t_i)(\mathbf v_0) + \gamma(\mathbf t_i)(\mathbf a) \mod U
  \]
 is in the transversal.
As $\cN$ is periodic with periodicity witnessed by the lattice $\bbbL$, $[ \mathbf v_1, \mathbf v_2 ] \in E$.
Let $\mathbf t_1'$ and $\mathbf t_2'$ be itineraries such that $\gamma(\mathbf t_i')(\mathbf v_0) = \mathbf v_i$ for $i = 1, 2$, 
 and these itineraries on $\cN'$ produce $\mathbf v_i'$ for $i = 1, 2$.
As $\{ \mathbf v_1, \mathbf v_2 \} \in E$, $[ \mathbf v_1', \mathbf v_2' ] \in E'$.
Choose $n_{i,1}, \ldots, n_{i,d} \in \bbbZ$  such that
  \[
 \gamma(\mathbf t_i)(\mathbf a) = \gamma(\mathbf t_i)(\mathbf a) \mod U = \sum_{j=1}^d n_{i,j} \mathbf l_j
  \]
 for $i = 1, 2, \ldots$, and the vertices
  \[
 \mathbf w_i = \mathbf v_i' + \sum_{j=1}^n n_{i,j} \mathbf l_j'
  \]
 satisfy $\{ \mathbf w_1, \mathbf w_2 \} \in E'$ by the periodicity of $\cN'$, so $\eta$ preserves edges.

The second claim is a sort of converse of the first.
Let $U'$ be the parallelopiped bounded by $\mathbf l_1', \ldots, \mathbf l_d'$:  $U'$ is the primitive unit cell 
 corresponding to $U$ for $\cN'$.
Suppose that $[\mathbf w_1, \mathbf w_2 ] \in E'$;  by periodicity, we get a corresponding $[ \mathbf v_1', \mathbf v_2' ] 
 \in E'$ where $\mathbf v_1'$ is within $U'$.
But the edge of $[ \mathbf v_1, \mathbf v_2 ] \in E$ is induced by the corresponding edge of $\cN^{\dagger}$ within U, and 
 the periodicity of $\cN$ does the rest.

Finally, to show that $\bbbS$ can be embedded in Sym$(\cN^{\dagger})$, let
  \[
 \bbbS' = \{ \1 \mathbf b', {}^{F^{-1}}M \2 : \1 \mathbf b, {}^{F^{-1}}M \2 \in \bbbS ]\},
  \]
 and we verify that $\bbbS'$ consists of symmetries of $\cN'$.
Suppose that $\1 \mathbf b, M \2 (\sum_{i=1}^d c_i \mathbf l_i + \mathbf v) = \sum_{i=1}^d c_i' \mathbf l_i + \mathbf w$ 
 where $\1 \mathbf b, M \2 (\mathbf v) = \sum_{i=1}^d c_i'' \mathbf l_i + \mathbf w$.
Then $M\sum_{i=1}^d c_i \mathbf l_i = \sum_{i=1}^d (c_i' - c_i'')\mathbf l_i$.
Using $\eta(\mathbf v) = \mathbf v'$ (for $\mathbf v \in U)$ and $\eta(\mathbf l_i) = \mathbf l_i'$ for $i = 1, \ldots, d$, 
 we have $\1 \mathbf b', M \2(\mathbf v') = \sum_{i=1}^d c_i'' \mathbf l_i' + \mathbf w'$ and:
\begin{eqnarray*}
 \1 \mathbf b', M \2 (\sum_{i=1}^d c_i \mathbf l_i' + \mathbf v') 
 &=& \mathbf b' + M\sum_{i=1}^d c_i \mathbf l_i' + M\mathbf v' \\
 &=& M\sum_{i=1}^d c_i \mathbf l_i' + \mathbf b' + M\mathbf v' \\
 &=& \sum_{i=1}^d (c_i' - c_i'') \mathbf l_i' + \1 \mathbf b, M \2 (\mathbf v') \\
 &=& \sum_{i=1}^d c_i' \mathbf l_i' + \mathbf w'.
\end{eqnarray*}
$\blacksquare$

\subsection{Deleting the Bad Nets}\label{degeneracy}

Following Definition~\ref{degenerate}, a solution would be ``degenerate'' if two edges of the original net $\cN$ 
 intersected.
Notice that degeneracy can be captured by systems of linear equations, as follows.
Formally, the edge $\{ \mathbf x, \mathbf y \}$ would intersect the edge $\{ \mathbf z, \mathbf w \}$ if there existed 
 $p, q \in [0, 1]$ such that 
  \[
 \mathbf x + p(\mathbf y - \mathbf x) = \mathbf z + q(\mathbf w - \mathbf z), 
  \]
 i.e., if
  \[
 (1 - p)\mathbf x + p\mathbf y - (1 - q)\mathbf z - q\mathbf w = \mathbf 0.
  \]
(Notice that if $p, q \in \{ 0, 1\}$, then this would characterize a collision of two vertices.)
For each $p, q \in [0, 1]$, this gives us a vector space to avoid, the result being that we want to take the nets 
 generated in Subsection~\ref{naive} and delete these vector subspaces.

More generally, notice that $\bigcup Q$ is a prefix-closed set of itineraries, and consider the following.
For each itinerary $\mathbf t \in \bigcup Q$ that is not maximal in $\bigcup Q$, and for each $(t, t')$ such that 
 $\mathbf t (t, t') \in \bigcup Q$, the last symbol $(t, t')$ represents an edge for the traveller to traverse.
Represent that edge $e$ with $d$-tuples of real variables $\mathbf x_{\mathbf t}$ and $\mathbf x_{\mathbf t(t, t')}$ 
 -- which we denote $\mathbf x_e$ and $\mathbf y_e$, respectively -- representing the endpoints of $e$.
Then, for each pair of distinct edges $e_1$, $e_2$ so represented in $\bigcup Q$, and for each pair of constants 
 $p_{e_1}, p_{e_2} \in [0, 1]$, we have the equation
  \begin{eqnarray}\label{ecollide}
 (1 - p_{e_1})\mathbf x_{e_1} + p_{e_1}\mathbf y_{e_1} - (1 - p_{e_2})\mathbf x_{e_2} - p_{e_2}\mathbf y_{e_2} = \mathbf 0.
  \end{eqnarray}
That includes, as special cases, for $\mathbf t_1, \mathbf t_2 \in \bigcup Q$ with $|\mathbf t_1|_{\simeq} \neq 
 |\mathbf t_2|_{\simeq}$, the equation
  \begin{eqnarray*}\label{vcollide}
 \mathbf x_{\mathbf t_1} - \mathbf x_{\mathbf t_2} = \mathbf 0.
  \end{eqnarray*}
Let $\cE^*$ be the set of all these equations generated from intersections of edges or collisions of vertices.

\begin{prop}\label{main2}
If $\cN'$ is a solution of $\cE$ but not a solution of any of the equations of $\cE^*$, then $\cN' \cong \cN$.
\end{prop}

\noindent{\bf Proof.}
By Proposition~\ref{main1}, it suffices to observe that the failure to satisfy any of the Equations~\ref{vcollide} 
 forces the homomorphism induced by the (equivalence classes of) itineraries to be one-to-one.
$\blacksquare$

\medskip

We would prefer that the edges don't intersect, which is guaranteed by the failure to satisfy any of the equations 
 of the Equations~\ref{ecollide}.

There are uncountably many vector subspaces to delete, so we should be a little careful.
We will use a straightforward observation from analysis.

\begin{prop}\label{closure}
Let $C \subseteq \bbbR^n$ be compact and $f {:} \; \bbbR^m \times \bbbR^n \to \bbbR$ be continuous.
For each $\mathbf c \in \bbbR^n$, let $\bbbW_{\mathbf c} = \{ \mathbf x {:} \; f(\mathbf x, \mathbf c) = 0 \}$.
Then $\bbbW = \bigcup_{\mathbf c \in C} \bbbW_{\mathbf c}$ is closed.
\end{prop}

\noindent{\bf Proof.}
It suffices to verify that $\bbbR^m - \bbbW$ is open.

Suppose that $\mathbf x \not\in \bbbW$;  we claim that $\inf \{ \| f(\mathbf x, \mathbf c) \| {:} 
 \; \mathbf c \in C \} > 0$.
Towards contradiction, suppose that there exists in $C$ a sequence $\mathbf c_j$, $j \to +\infty$, with 
 $\|f(\mathbf x, \mathbf c_j)\| \to 0$.
By the compactness of $C$, there is in $C$ an accumulation point $\mathbf c_{\infty}$ of the sequence $\mathbf c_j$, 
 $j \to +\infty$, so by the continuity of $f$ we get $f(\mathbf x, \mathbf c_{\infty}) = 0$ and $\mathbf x \in \bbbW$ 
 after all.

So if $\mathbf x \not\in \bbbW$ and we set $\varepsilon = \inf\{\|f(\mathbf x, \mathbf c)\| {:} \; \mathbf c \in C \}
 > 0$, by the continuity of $f$ we can choose $\delta > 0$ so that for all $\mathbf y \in \bbbR^m$,
  \[
 \| \mathbf x - \mathbf y \| < \delta \implies \|f(\mathbf x, \mathbf c) - f(\mathbf y, \mathbf c) \| < \varepsilon.
  \]
Thus for all such $\mathbf y$, $\mathbf y \not\in\bbbW$.
So all points $\mathbf x$ in $\bbbR^m - \bbbW$ are interior points, so $\bbbR^m - \bbbW$ is open, so $\bbbW$ is closed.
$\blacksquare$

\medskip

We pull all this together.

\begin{theo}
For each periodic net, there is an isomorphic periodic net whose vertices are integer points (modulo an appropriate 
 affine transform).
Thus $\nu$ maps the symmetry group $\bbbS$ of $\cN$ to a subgroup of the symmetry group $\bbbS'$ of $\cN'$, and hence every 
 orbit of $\cN$ is a (not necessarily proper) suborbit of an orbit of $\cN'$.
\end{theo}

Thus the symmetry group $\bbbS$ of $\cN^{\dagger}$ is isomorphic to a (not necessarily proper) subgroup of $\bbbS'$ of $\cN'$, so 
 in particular the point group of $\cN^{\dagger}$ is a (not necessarily proper) subgroup of the point group of $\cN'$.

\medskip

\noindent{\bf Proof.}
Let $\cN^{\dagger}$ be a scaffolded net obtained from a periodic net $\cN$, and let its point group be a subgroup of the 
 maximal point group $\bbbH$, conjugate to the integer matrix group ${}^{F^{-1}}\bbbH^*$.
By Proposition~\ref{DD}, $\cN$ admits a fundamental transversal and using it and matrices from $\bbbH^*$, we can generate 
 a quotient graph $Q$ as in Subsection~\ref{naive}.

From the quotient graph $\cN^{\dagger}$, we get a system of equations $\cE$ whose solutions give the positions of vertices 
 and edges of a unit cell that are homomorphic images of $\cN^{\dagger}$;  
 this set of solutions forms a vector space $\bbbE$.
Meanwhile, we also get a system of equations $\cE^*$ such that any solution of $\cE$ represented by a tuple of vertices 
 in $\bbbE$ that satisfies any equation of $\cE^*$ must be degenerate;  conversely, any solution of 
 $\cE$ that does not satisfy any equation of $\cE^*$ is not degenerate and hence by Proposition~\ref{main2} 
 must represent a net isomorphic to $\cN$, along with symmetries corresponding to $\Sym(\cN)$.

Let $\bbbE^*$ be the set of all tuples of $\bbbE$ that satisfy some equation of $\cE^*$, and notice that $\cE^*$ consists 
 of homogeneous linear equations that are of one of two forms.

The kind of equation of $\cE^*$, corresponding to vertex collisions, is
  \[
 \mathbf x_i - \mathbf x_j = \mathbf 0,
  \]
 and there are only finitely many of these, and each of their solution spaces $\bbbW_{i,j}$ is a vector space, hence 
 closed.
It follows that the finite union 
  \[
 \bbbW' = \bbbW \cup \bigcup_{i,j {:} \; \mbox{\scriptsize ``$\mathbf x_i - \mathbf x_j = \mathbf 0$''} \in \cE^*} \bbbW_{i,j}
  \]
 is closed.

The edge collisions are represented by equations of the form
  \[
 \sum_i c_i \mathbf x_i + \sum_j c_j' \mathbf y_j = \mathbf 0
  \]
 where the set of all tuples $\mathbf c = (c_1, \ldots, c_1', \ldots)$ satisfies $c_i, c_j' \geq 0$ for all $i$, $j$, 
 and where there exist $i_1, i_2, j_1, j_2$ such that
  \[
 i \not\in \{ i_1, i_2 \} \implies c_i = 0 \quad \& \quad j \not\in \{ j_1, j_2 \} \implies c_j = 0
  \]
 and $c_{i_1} + c_{i_2} = c_{j_1} + c_{j_2} = 1$.
This set of tuples $\mathbf c$ is compact, so by Proposition~\ref{closure}, the union $\bbbW$ of all the solution spaces 
 $\bbbW_{\mathbf c}$ of equations of the first form is closed.

Thus the space $\bbbE - \bbbW'$ of solutions representing nets isomorphic to $\cN$ is open.

As the tuple of vertices of $\cN^{\dagger}$'s unit cell is an element of $\bbbE - \bbbW'$, $\bbbE - \bbbW' 
 \neq \varnothing$.
As $\bbbE - \bbbW'$ is open, there exists a neighborhood $N$ around $\cN^{\dagger}$'s tuple wholly within $\bbbE - \bbbW'$.
As all the coefficients of equations of $\cE$ are integers, $\bbbE$ is spanned by integer vectors, and hence the set of 
 tuples of tuples of rationals is dense in $\bbbE$.
From Observation~\ref{coeffs}, there is a tuple of tuples of rationals in $N$, i.e., there exists a net $\cN'$ whose vertices 
 in its unit cell, and making up its lattice, are all rational;  thus by periodicity, all the vertices of $\cN'$ are rational.
Expressing all coordinates of vertices and lattice vectors of $\cN'$ as reduced fractions, there are finitely many 
 integers appearing in the denominators, so we can choose the least common multiple $m$ of all denominators appearing 
 in these fractions, and multiply every vertex vector by $m$ to get a net $m \cN'$ isomorphic to $\cN$ and whose 
 vertices are all at integer points.

Now we make the final change of basis back get the desired scaffolded net $\cN^{\dagger}$.
$\blacksquare$

\medskip

Having obtained the net $\cN'$ of integer points, isomorphic to $\cN$ and whose symmetry group admits a subgroup 
 isomorphic to $\Sym(\cN)$, we have completed the proof of Theorem~\ref{main}.

\section{Conclusion}\label{coda}

The original motivation for this paper was the development of a ``Crystal Turtlebug'' algorithm that enumerates 
 crsytal graphs with vertices on a geometric lattice.
The question was whether this algorithm eventually enumerates representatives (of maximal symmetry) of all 
 isomorphism classes of crystal graphs.
One consequence of Theorem~\ref{main} is that it does.

The Crystal Turtlebug is a project in a growing field of mathematical and computational crystallography.
Although graphical representations of molecules and solids go back to the Nineteenth century, it is only 
 in the last few decades that a systematic attempt has been made to design, organize, catalogue and apply 
 graphical representations of specific crystals.
The thread of research arising from \cite{Wells} through \cite{OH} to recent works for chemists like 
 \cite{Ohrstrom} and \cite{Lord} is motivated by a desire to design crystals in advance prior to synthesis, 
 rather than relying on combinatorial chemistry to physically search through thousands of initial conditions 
 in the hope of finding one that it interesting (\cite{ferey}, \cite{yaghi}).

Recently, several groups have composed computer programs that generate crystal nets in the hope that they 
 may prove to be viable blueprints.
Some groups have developed algorithms based on fundamentally geometric principles, e.g., by enumerating tilings of 3-space 
 (\cite{DDHKM99}), by reflecting a fundamental region around (\cite{TRBRF04}), or by attaching polyhedra together one at 
 a time (\cite{LeBail}).
(There are groups employing more distant algorithms, e.g., \cite{Hyde} and \cite{Deem};  see \cite{survey1} and 
  \cite{survey2} for more.)
And many of the catalogues (e.g., the library attached to SYSTRE (of the Generation, Analysis and Visualization of Reticular 
 Ornaments using GAVROG \cite{SYSTRE}), the Reticular Chemistry Structure 
 Resource (\cite{RCSR}), and TOPOS (see, e.g., \cite{TOPOS} or go to http://www.topos.ssu.samara.ru/)) use crystal net 
 isomorphism as a principle standard of identification, the issue of whether a crystal net is novel depends on whether 
 it is isomorphic to any extant crystal nets.

For crystal design, then, the message of this paper is that if all one desires is to generate isomorphism classes of crystal 
 nets, it suffices to generate (representative) nets with integer points as vertices.
In addition, because we are usually interested in the symmetry groups of these nets (i.e., isometry groups that induce 
 automorphisms on these nets), for a given isomorphism class of crystal nets, we can generate such a representative net so 
 that its point group is maximal among the point groups of nets within this isomorphism class.

This result is not surprising, considering Bieberbach's ``Second'' Theorem (\cite{Bieberbach1}, see \cite{charlap}) that 
 every crystallographic group is an affine conjugate of a group whose orbit of the origin consists of integer 
 points.
Indeed, the main result of this paper is a generalization of the Second Theorem:
\begin{itemize}
 \item
  Given any crystallographic group $\bbbG$ whose subgroup of translations is generated by translation of vectors $\mathbf l_1, 
   \cdots, \mathbf l_n$, start with the transversal consisting of a vertex at the origin and edges from the origin to 
   adjacent vertices (of the same orbit) at $\mathbf l_1, \cdots, \mathbf l_n$.
 \item
  By Theorem~\ref{main}, there is another Euclidean graph isomorphic to the one just constructed, whose symmetry group 
   is represented by matrices of integers (possibly modulo an affine transformation) and having a subgroup $\bbbG'$ 
   isomorphic to $\bbbG$.
\end{itemize}
Then $\bbbG'$ is the desired conjugate.
 
The intended application of this theorem was to verify that if a computer program will generate a (unit cell of) any 
 crystal net whose vertices are integer points (modulo the appropriate affine transformation\footnote{For physical 
 crystallography, the transformations are either the identity (if the crystal's point group is a subgroup of $m\bar 3m$) 
 or the map generated by the assignment of basis elements $\langle 1, 0, 0\rangle \mapsto \langle 1, 0, 0\rangle$, 
 $\langle 0, 1, 0 \rangle \mapsto \langle 1/2, \sqrt 3/2, 0 \rangle$, $\langle 0, 0, 1 \rangle \mapsto \langle 0, 0, 1 
 \rangle$ (if the crystal's point group is a subgroup of $6/mmm$).} if necessary).

This project started with a computing project proposed by the author to W. E. Clark in 2007, who composed 
 a sequence of programs in MAPLE, one of which was generating novel uninodal nets by early 2008, and 
 whose underlying rationale is explained in \cite{Clark}, resting on the results on vertex transitivity in 
 \cite{sabidussi}.
The author modified this algorithm to obtain a program for binodal edge transitive nets, and conceptually for 
 any net (\cite{mainpaper}, the latter program being under development).

The author is grateful in particular to W. E. Clark for his assistance and advice, and to the University of 
 South Florida for providing a sabbatical during the spring of 2008, during which the author composed the 
 first versions of the program that is now generating nets for chemists to try to realize in the lab.

\end{document}